\documentclass[12pt]{amsart}  
\usepackage{amscd} 
\usepackage{mathrsfs}
\newtheorem{theorem}{Theorem}[section]

\theoremstyle{definition}

\newtheorem{remark}[theorem]{Remark} 
\numberwithin{equation}{section}

\def\H{\mathscr H}  
\def\K{\mathscr K}

\def\D{\mathscr D} 
\def\R{\mathscr R}
\def\Ex{\mathscr E}
\def\RE{\mathbb R} 
\def\CO{{\mathbb C}}
\def\NA{\mathbb N} 
\def\IN{{\mathbb Z}}

\def\uno{\mathsf 1}

\def\fh{\mathfrak h}

\def\E{\mathsf {E}}
\def\P{{\mathsf P}}

\def\min{\text{\rm min}}
\def\max{\text{\rm max}}
\newsymbol\restriction 1316
\newsymbol\subsetneq 2328
                                                                    
\begin{document}

\title[Kre\u\i n's
Formula for Elliptic Differential Operators]
{Krein's 
Resolvent Formula for Self-Adjoint Extensions of 
Symmetric Second Order Elliptic Differential Operators}

\author{Andrea Posilicano}
\author{Luca Raimondi}
\address{Dipartimento di Scienze Fisiche e Matematiche,  Universit\`a
dell'Insubria, I-22100 Como, Italy}

\email{posilicano@uninsubria.it}
\email{luca.raimondi@yahoo.it}

\thanks{
(primary), 
47B38, 35J25 (secondary)}

\begin{abstract}
Given a symmetric, semi-bounded, second order elliptic differential
operator $A$ on a bounded domain with $C^{1,1}$ boundary, we provide a
Kre\u\i n-type formula for the resolvent difference between its
Friedrichs extension and an arbitrary self-adjoint one. 
\end{abstract}
{\maketitle }
\begin{section}{Introduction.}

Given a bounded open set $\Omega\subset\RE^n$, $n>1$,  let us consider
a second order elliptic differential operator
$$
A:C^\infty_c(\Omega)\subset L^2(\Omega)\to  L^2(\Omega)\,,\quad
A=\sum_{i,j=1}^n \partial_i(a_{ij}\partial_{j})-
\sum_{i=1}^nb_i\partial_i-c\,.
$$ Such an operator $A$, under appropriate hypotheses on its
coefficients and on $\Omega$ (these will be made precise in the
section 3), is closable and its closure $A_{\min}$, the minimal
realization of $A$, has domain given by $H^2_0(\Omega)$, the closure
of $C^{\infty}_c(\Omega)$ with respect to $H^2(\Omega)$ Sobolev norm.
If $A$ is symmetric then $A_{min}$ is symmetric but not
self-adjoint, i.e.  $A$ is not essentially self-adjoint.  Indeed
$A_\min^*=A_\max$, where $A_\max$, the maximal realization of $A$, has
domain made by the functions $u\in L^2(\Omega)$ such that $Au\in
L^2(\Omega)$.  Assuming that $A_\min$ is semibounded, then $A_\min$
has a self-adjoint extension $A_0$ (the Friedrichs extension,
corresponding to Dirichlet boundary conditions), $A_\min\subsetneq
A_0\subsetneq A_\max$, and hence $A_\min$ has infinitely many self-adjoint extensions.  \par The problem of the parametrization of all
self-adjoint extensions of $A_\min$ in terms of boundary conditions
was completely solved (in the case of an elliptic differential
operator of arbitrary order) in \cite{[G]} (for some older papers about similar topics we just quote \cite{[C]} and \cite{[V]}).  Here, by using the
approach developed in \cite{[Po1]}-\cite{[Po4]}, we give an
alternative derivation of such a result by providing a Kre\u\i n-like
formula for the resolvent difference between an arbitrary self-adjoint
extension of $A_\min$ and its Friedrichs extension $A_0$. For the sake
of simplicity here we consider the case of a second order differential
operator. The case of higher order operators can be treated in a
similar way.  \par 
In the case $A$ is the Laplacian, the Kre\u\i n resolvent formula here presented 
has been given in \cite{[Po4]}, Example 5.5.  For other recent
results on Kre\u\i n-type formula for partial differential operators
see \cite{[BL]}, \cite{[Ry]}, \cite{[BMNW]}, \cite{[GM1]},
\cite{[GM2]}, \cite{[BGW]}.  \par 
In order to help the reader's intuition on the results here presented, in Section 4 we consider one of the simplest possible examples: a rotation-invariant elliptic operators $A$ on the disc $D\subset
\RE^2$. Thus, notwithstanding the symmetric operator here considered has infinite deficiency indices, 
due to the presence of symmetries the resolvents of their self-adjoint extensions can be written, 
by separation of variables, in a form which resembles the finite indices case (see the comments in Remark 4.1), and the 
corresponding spectral analysis becomes simpler. 
As illustration, given any sequence $\{\lambda_n\}_1^\infty\subset \RE$, boundary conditions at $\partial D$ can be given for which
$A$ is self-adjoint and such that $\{\lambda_n\}_1^\infty$ is
contained in its point spectrum. Remark 
\ref{remark} shows that such
boundary conditions can be quite different from the usual ones.

\end{section}

\begin{section}{Preliminaires}
For the reader's convenience in this section we collect some
results from \cite{[Po1]}-\cite{[Po4]}. We refer to these papers, in particular to \cite {[Po4]}, for a through discussion about the 
connection of the approach here presented with both the standard von Neumann's theory of self-adjoint extension \cite{[vN]} and with Boundary Triple Theory \cite{[DM]},\cite{[GG]}. 
\par
From now on we will denote by
$$\D(L)\,,\quad\K(L)\,,\quad\R(L)\,,\quad\rho(L)$$ the domain, kernel,
range and resolvent set of a  linear operator $L$.  \par Let
$\mathscr{H}$ be a Hilbert space with scalar product
$\langle\cdot,\cdot\rangle$ and let  $$
A_0:\D(A_0)\subseteq\mathscr{H}\rightarrow\mathscr{H}
$$
a self-adjoint operator on it. We denote by  $\mathscr{H}_{A_0}$ be
the Hilbert space given by the linear space $\D(A_0)$ endowed with the
scalar product
$$
\langle\phi,\psi\rangle_{A_0}=\langle\phi,\psi\rangle+\langle
A_0\phi,A_0\psi\rangle\,.
$$
Given then a Hilbert space $\mathfrak{h}$ with scalar product
$(\cdot,\cdot)$ and a linear, bounded and surjective operator
$$\tau:\H_{A_0}\rightarrow\mathfrak{h}\,,$$ such that $\K(\tau)$ is
dense in $\H$, we denote by $S$ the densely defined closed symmetric
operator
$$
S:\K(\tau)\subseteq\H\to\H\,,\qquad S\phi:=A_0\phi\,.
$$
Our aim is to provide, together with their resolvents, all
self-adjoint extension of $S$. \par For any $z\in\rho(A_0)$ we define
the bounded operators
\begin{equation*}
R_z:=(-A_0+z)^{-1}:\H\to\H_{A_0}\,,
\end{equation*}
\begin{equation}\label{gz}
G_z:=(\tau R_{\bar z})^*:\mathfrak{h}\rightarrow \H\,.
\end{equation}
By \cite{[Po2]}, Lemma 2.1, given the surjectivity hypothesis
$\R(\tau)=\fh$, the density assumption $\overline{\K(\tau)}=\H$ is equivalent
to
\begin{equation*}
\R(G_z)\cap \D(A_0)=\{0\}\,.
\end{equation*}
However, since by first resolvent identity
\begin{equation}\label{GR1}
(z-w)R_wG_z=G_w-G_z\,,\\
\end{equation}
one has
\begin{equation*}
\R(G_w-G_z)\subset \D(A_0)\,.\\
\end{equation*}
From now on, even if this hypothesis can be avoided (see
\cite{[Po1]}-\cite{[Po4]}),  for the sake of simplicity we suppose that
$$
0\in\rho(A_0)\,.
$$ 
We define the family $\Gamma_z$, $z\in\rho(A_0)$, of bounded linear
maps
\begin{equation}\label{gamma}
\Gamma_z:\fh\rightarrow\fh\,,\qquad
\Gamma_z:=\tau\left(G_0-G_z\right)\equiv -z\tau A_0^{-1}G_z \,.
\end{equation}
Given then an orthogonal projection
$$
\Pi:\mathfrak{h}\rightarrow \fh\,,\qquad \fh_0\equiv\R(\Pi)\,,
$$
and a self-adjoint operator
$$
\Theta:\D(\Theta)\subseteq \mathfrak{h}_0\rightarrow \mathfrak{h}_0,
$$
we define the closed operator
\begin{equation*}
\Gamma_{z,\Pi,\Theta}:=(\Theta+\Pi\,\Gamma_z\,\Pi):\D(\Theta)\subseteq
\mathfrak{h}_0\rightarrow \mathfrak{h}_0,
\end{equation*}
and the open set
$$
Z_{\Pi,\Theta}:=\{z\in\rho(A_0)\::\:0\in\rho(\Gamma_{z,\Pi,\Theta})\}.
$$
With such premises the next two theorems have  straightforward
proofs. Theorem \ref{aggiunto} is an obvious modification (taking
into account the hypothesis  $0\in\rho(A_0)$) of  Theorem 3.1 in
\cite{[Po3]} (also see  \cite{[Po2]}, Theorem 3.4);  Theorem
\ref{estensioni} is the  combination of Theorem  \ref{aggiunto} with
Theorem 2.1 and Theorem 2.4 in \cite{[Po4]}  (also see \cite{[Po1]},
Theorem 2.1, \cite{[Po2]}, Theorem 2.2, for the case $\Pi=\uno$).
\begin{theorem}\label{aggiunto} The adjoint of $S$ is given by
$$
S^*:\D(S^*)\subseteq\H\to\H\,,\qquad S^*\phi=A_0\phi_0\,,
$$
$$
\D(S^*)=\{\phi\in\H\,:\,\phi=\phi_0+G_0\zeta_{\phi},\  \phi_0\in
\D(A_0),\ \zeta_{\phi}\in\mathfrak{h}\}\,.
$$
Moreover
\begin{equation}\label{Green}
\forall\,\phi\,,\psi\in\D(S^*)\,,\quad\langle S^*\phi,\psi\rangle-\langle\phi,S^*\psi\rangle =
(\tau\phi_0,\zeta_{\psi})-(\zeta_{\phi},\tau\psi_0)\,.
\end{equation}
\end{theorem}

\begin{theorem}\label{estensioni}
The set $Z_{\Pi,\Theta}$ is not void,
$$\mathbb{C}\setminus \mathbb{R}\subseteq Z_{\Pi,\Theta}$$ and
\begin{equation*} R_{z,\Pi,\Theta}
:=R_z+G_z\Pi\,\Gamma_{z,\Pi,\Theta}^{-1}\,\Pi\,  G_{\bar z}^*\,,\qquad
z\in Z_{\Pi,\Theta}\,,
\end{equation*}
is the resolvent of the self-adjoint extension ${A}_{\Pi,\Theta}$ of
$S$ defined by
$$
A_{\Pi,\Theta}:\D(A_{\Pi,\Theta})\subseteq\H\to\H\,,\qquad
A_{\Pi,\Theta}\phi=S^*\phi\equiv A_0\phi_0\,,
$$
$$
\D(A_{\Pi,\Theta})= \{\phi\in\D(S^*)\,:\,\zeta_{\phi}\in\D(\Theta)\,,\
\Pi\tau\phi_0=\Theta\zeta_\phi\}\,.
$$
\end{theorem}
\begin{remark}\label{regtau}
Note that, since $\phi_0=A_0^{-1}S^*\phi$,
$$\Pi\tau\phi_0=\Theta \zeta_\phi\quad\iff\quad
\Pi\hat\tau_0\phi=\Theta \zeta_\phi\,,
$$
where the regularized trace operator $\hat\tau_0$ is defined  by
$$
\hat\tau_0:\D(S^*)\to\fh\,,\quad\hat\tau_0\phi :=\tau
A_0^{-1}S^*\phi\,.
$$
\end{remark}
By exploiting the connection with  von Neumann's theory (see
\cite{[Po4]}, section 3; see also \cite{[Po2]},  section 4, for the
case of relatively prime extensions) one obtains
\begin{theorem}\label{bundle}
The set of operators provided by Theorem \ref{estensioni}  coincides
with the set $\Ex(S)$ of all self-adjoint extensions of the symmetric
operator $S$.  Thus $\Ex(S)$ is parametrised by the bundle
$p:\E(\fh)\to\P(\fh)$, where $\P(\fh)$ denotes the set of orthogonal
projections in $\fh$ and $p^{-1}(\Pi)$ is the set of self-adjoint
operators in the Hilbert space $\R(\Pi)$.  The set of self-adjoint operators in
$\fh$, i.e. $p^{-1}(\uno)$,  parametrises all  relatively prime
extensions of $S$ i.e. those for which $\D(\hat A)\cap
\D(A_0)=\D(S)$.
\end{theorem}
We conclude this section with a result about the spectral properties
of the extensions (see \cite{[DM]}, Section 2, for point 1 and
\cite{[Po3]}, Theorem 3.4, for point 2).
\begin{theorem}\label{spettro}
1)$$ \lambda\in\sigma_p(A_{\Pi,\Theta})\cap\rho(A_0)\quad\iff\quad
0\in\sigma_p(\Gamma_{\lambda,\Pi,\Theta}),
$$
where $\sigma_p(\cdot)$ denotes point spectrum. An analogous result
holds for the continuous spectrum.  \par\noindent  2)
$$
G_\lambda:\K(\Gamma_{\lambda,\Pi,\Theta})\rightarrow
\K(-A_{\Pi,\Theta}+\lambda)
$$
is a bijection for any
$\lambda\in\sigma_p(A_{\Pi,\Theta})\cap\rho(A_0)$.
\end{theorem}

\end{section}

\begin{section}{Extensions and Krein's Formula.}
Let $\Omega\subset\RE^n$, $n>1$,  a bounded open set with a Lipschitz
boundary. We denote by $H^k(\Omega)$ the Sobolev-Hilbert  space given
by closure of  $C^{\infty}(\bar\Omega)$ with respect to the norm
$$
\|u\|_{H^k(\Omega)}^2=\sum_{0\le \alpha_1+\dots+\alpha_n\le k} \|
\partial_1^{\alpha_1}\dots \partial_n^{\alpha_n} u\|^2_{L^2(\Omega)}\,.
$$
Analogously $H^k_0(\Omega)\subsetneq H^k(\Omega)$ denotes  the closure
of $C_c^{\infty}(\Omega)$ with respect to the  same norm.\par Given
the differential expression
$$
A=\nabla\!\cdot\!a\nabla -b\!\cdot\!\nabla-c\equiv \sum_{i,j=1}^n
\partial_i(a_{ij}\partial_{j})- \sum_{i=1}^nb_i\partial_i-c
$$
we suppose that  the matrix $a(x)\equiv(a_{ij}(x))$ is Hermitean for
a.e.  $x\in\Omega$, that there exist $\mu_1>0$, $\mu_2>0$  such that
$$\forall\xi\in\RE^n\,,\qquad \mu_1\|\xi\|^2\le \xi\!\cdot\!a(x)\xi
\le \mu_2\|\xi\|^2$$  and that
$$b_i\in L^{q}(\Omega)\,,\quad c\in L^{q/2}(\Omega)\,, \quad q=n
\ \text{\rm if}\ n\ge 3\,,\ q>2\ \text{\rm if}\ n=2\,.$$ Then $A$ maps
$H^1(\Omega)$ into $H^{-1}(\Omega)$ (see e.g. \cite{[EE]}, Section 1,
Chapter VI), where $H^{-1}(\Omega)$ denotes the adjoint space of
$H^{1}_0(\Omega)$, the sesquilinear form
$$ q_A:H^1_0(\Omega)\times H^1_0(\Omega)\to\CO
$$
$$ q_A(u,v):= -\left(\langle\nabla u,a\nabla v\rangle_{L^2(\Omega)}
+\langle u,b\!\cdot\!\nabla v\rangle_{L^2(\Omega)} +\langle
u,cv\rangle_{L^2(\Omega)}\right)\,,
$$ is continuous and there exists a positive constant $\lambda$ such
that $-q_A+\lambda$ is coercive (see e.g.  \cite{[EE]}, Proposition
1.2, Chapter VI).  Thus by Lax-Milgram Theorem (see e.g. \cite{[EE]},
Theorem 1.4, Chapter VI) there exists a unique closed, densely
defined, linear operator
$$ A_0:\D(A_0)\subseteq L^2(\Omega)\to L^2(\Omega)\,,\qquad A_0u=Au\,,
$$
$$ \D(A_0)=\{u\in H^1_0(\Omega)\,:\, Au\in L^2(\Omega)\}\,,
$$ such that
$$ \forall u\in \D(A_0)\,,\quad\forall v\in H^1_0(\Omega)\,,\qquad
q_A(u,v)=\langle u,A_0v\rangle_{L^2(\Omega)}\,.
$$ Moreover $\D(A_0)$ is dense in $H^1_0(\Omega)$, $0\in \rho
(-A_0+\lambda)$, $A_0$ has a compact resolvent and its spectrum
consists of an infinite sequence of eigenvalues
$\lambda_n$, each having finite multiplicity and with
Re$\lambda_n<-\lambda$.  An analogous result holds for the
sesquilinear form $q_A^*$,
$$q_A^*(u,v):=\overline {q_A(v,u)}$$ and the operator corresponding to
$q_A^*$ is the adjoint $A_0^*$. \par Suppose now that
$$\partial_ia_{ij}\in L^{q}(\Omega) \,, \quad q=n \ \text{\rm
  if}\ n\ge 3\,,\ q>2\ \text{\rm if}\ n=2\,,
$$ so that, by Sobolev Embedding Theorem, $A$ is continuous from
$H^2(\Omega)$ into $L^2(\Omega)$ and
$$H^2_0(\Omega)\subsetneq H^2(\Omega)\cap
H^1_0(\Omega)\subseteq\D(A_0)\,.$$ By interior regularity estimates
(see e.g.  \cite{[LU]}, Section 7, Chapter 3)
$A\!\!\restriction_{C_c^\infty(\Omega)}$, the restriction of $A$ to
$C_c^\infty(\Omega)$, is closable and its closure is given by
$A_{\min}\subsetneq A_0$, the minimal realization of $A$, defined by
$$ A_{\min}:H^2_0(\Omega)\subseteq L^2(\Omega)\to L^2(\Omega)\,,\qquad
A_{\min}u:=Au\,.
$$ From now on we suppose that $$q_A=q_A^*\,.$$ Thus $A_0$ is a
self-adjoint operator, the Friedrichs extension of the closed
symmetric operator $A_\min$ and one has
$$A_\min^*=(A\!\!\restriction_{C_c^\infty(\Omega)})^*=A_\max\,,$$
where $A_{\max}$, the maximal realization of $A$, is defined by
$$ A_{\max}:\D(A_\max)\subseteq L^2(\Omega)\to L^2(\Omega)\,,\qquad
A_{\max}u:=Au\,,
$$
$$ \D(A_\max):=\{u\in L^2(\Omega)\,:\, Au\in L^2(\Omega)\}\,.
$$ Hence $$\D(A_0)=H^1_0(\Omega)\cap\D(A_\max)\,.$$ Moreover
$$\D(A_\min)=H^2_0(\Omega)\subsetneq\D(A_\max)\,,$$ so that
$A\!\!\restriction_{C_c^\infty(\Omega)}$ is not essentially
self-adjoint, $$A_\min\subsetneq A_0\subsetneq A_\max\,,$$ and the
symmetric operator $A_\min$ has infinitely many self-adjoint extensions. We
want now to find all such extensions and to give their resolvents. In
order to render straightforward the application of the results given
in Section 2, we would like to have a more explicit characterization
of $\D(A_0)$.  Thus in the following we impose more stringent
hypotheses on the set $\Omega$.  \par Suppose that the boundary of
$\Omega$ is a piecewise $C^2$ surface with curvature bounded from
above and that $a_{ij}\in C(\bar\Omega)$ when $n\ge 3$. Then, by
global regularity results (see e.g. \cite{[LU]}, Chapter 3, Section
11), the graph norm of $A_\max$ is equivalent to that of 
$H^2(\Omega)$ on $C_0^\infty(\bar\Omega)$, the space of smooth functions on $\Omega$
which vanish on its boundary $\partial\Omega$.  Thus
$A\!\!\restriction_{C_0^\infty(\bar\Omega)}$, the restriction of $A$
to $C_0^\infty(\bar\Omega)$, is closable and its closure is given by
$$ \tilde A_{0}:\tilde H^2_{0}(\Omega)\subseteq L^2(\Omega)\to
L^2(\Omega)\,,\qquad \tilde A_{0}u:=Au\,,
$$ where $\tilde H^2_{0}(\Omega)$ denotes the closure of
$C^\infty_0(\bar\Omega)$ with respect to the $H^2(\Omega)$ norm. \par
Without further hypotheses on $\Omega$, $\tilde A_{0}\not=A_0$ is
possible: for example if $\Omega$ is a non-convex plane polygon then the Laplace
operator $\Delta$ is not self-adjoint on $\tilde
H^2_{0}(\Omega)$. Indeed by \cite{[BS]} it has deficiency indices $(d_-,d_+)=(d,d)$,
where $d$ is the number of non-convex corners .  \par Suppose now that
the $a_{ij}$'s are Lipschitz continuous up to the boundary and that
$\partial\Omega$ is $C^{1,1}$, i.e. it is locally the graph of a $C^1$
function with Lipschitz derivatives (see e.g \cite{[Gr]}, Section 1.2,
for the precise definition). Then (see e.g. \cite{[LM]}, Chapter 1,
Section 8.2, \cite{[Gr]}, Section 1.5) there are unique continuous and
surjective linear maps
$$ \rho:H^1(\Omega)\to H^{1/2}(\partial\Omega)\,,
$$
$$ \gamma_a:H^2(\Omega)\to H^{3/2}(\partial\Omega)\oplus
H^{1/2}(\partial\Omega)\,,\qquad
\gamma_a\phi:=\left(\rho\phi,\tau_a\phi\right)\,,
$$ such that
$$ \rho\phi\,(x):=\phi\,(x)\,,\quad\tau_a\phi\,(x)\equiv
\frac{\partial \phi}{\partial \nu_a}\,(x):=\sum_{i,j=1}^{n} a_{ij}(x)
\nu_i(x){}^{}\partial_j\phi\,(x)
$$ for any $\phi\in C^\infty(\bar\Omega)$ and $x\in\partial\Omega$.  Here
$\nu\equiv(\nu_1,\dots,\nu_n)$ denotes the outward normal vector on
$\partial\Omega$ and $H^s(\partial\Omega)$, $s>0$, are the usual
fractional Sobolev-Hilbert spaces on $\partial\Omega$ (see
e.g. \cite{[Gr]}, Section 1.3.3).  Moreover Green's formula holds: for
any $u\in H^2{(\Omega)}$ and $v\in H^2{(\Omega)}\cap H^1_0{(\Omega)}$
one has
\begin{equation}\label{Green1}
\langle Au,v\rangle_{L^2(\Omega)}= \langle u,A_0v\rangle_{L^2(\Omega)}
-                \langle\rho u,\tau_av\rangle_{L^2(\partial\Omega)}\,.
\end{equation}
By proceeding as in the proof of Theorem 6.5 in \cite{[LM]}, Chapter 6
(which uses (\ref{Green1})) the map $\gamma_a$ can be extended to (see
\cite{[Gr]}, Theorem 1.5.3.4)
$$ \hat\gamma_a:\D(A_{\max})\to H^{-1/2}(\partial\Omega)\oplus
H^{-3/2}(\partial\Omega)\,,\qquad \hat \gamma_a\phi
=(\hat\rho\phi,\hat\tau_a\phi)\,,
$$ where $H^{-s}(\partial\Omega)$ denotes the adjoint space of
$H^{s}(\partial\Omega)$, and Green's formula (\ref{Green1}) can be
extended to the case in which $u\in \D(A_\max)$:
\begin{equation}\label{Green2}
\langle A_\max u,v\rangle_{L^2(\Omega)}= \langle
u,A_0v\rangle_{L^2(\Omega)} -(\hat\rho
u,\tau_av)_{-\frac{1}{2},\frac{1}{2}}\,.
\end{equation}
Here $(\cdot,\cdot)_{-\frac{1}{2},\frac{1}{2}}$ denotes the duality
between $H^{1/2}(\partial\Omega)$ and $H^{-1/2}(\partial\Omega)$.
With such definitions of $\rho$ and $\tau$ one has (see
e.g. \cite{[Gr]}, Corollary 1.5.1.6),
$$ H^1_0(\Omega)=H^1(\Omega)\cap\K(\rho)\,,\qquad
H^2_{0}(\Omega)=H^2(\Omega)\cap\K(\gamma_\uno)\,.
$$ Moreover, by the stated properties of $\rho$ and $\hat\rho$, by the
equivalence of the graph norm of $A_\max$ with the $H^2(\Omega)$ norm
on $\tilde H^2_0(\Omega)$ and by the density of $C^\infty(\bar\Omega)$
in $\D(A_\max)$, one gets the equalities
$$\tilde H^2_{0}(\Omega)=H^2(\Omega)\cap H^1_0(\Omega)= \D(A_\max)\cap
H^1_0(\Omega)\equiv \D(A_0)\,,$$ so that $\tilde A_{0}=A_0$. \par In
conclusion we can apply the results given in Section 2 (by adding, if necessary, a constant to $A_0$ we may suppose that $0\in\rho(A_0)$) to the
self-adjoint operator
$$ A_0:H^2(\Omega)\cap H^1_0(\Omega)\subseteq L^2(\Omega)\to
L^2(\Omega)\,,\qquad A_0u:=Au\,,
$$ with $S=A_\min$, $\fh=H^{1/2}(\partial\Omega)$ and
$$ \tau:H^2(\Omega)\cap H^1_0(\Omega)\to
H^{1/2}(\partial\Omega)\,,\qquad \tau:=\tau_a
$$ Note that $\K(\tau)=H^2_0(\Omega)$ since
$\K(\gamma_a)=\K(\gamma_\uno)$ by $\nu(x)\!\!\cdot\!a(x)\nu(x)\ge
\mu_1>0$, $x\in\partial\Omega$, and that $\tau$ is surjective by the
surjectivity of $\gamma_a$.  \par Thus, by Theorem \ref{bundle}, under
the hypotheses above, the set $\Ex(A_\min)$ of all self-adjoint
extensions of $A_\min$ can be parametrized by the bundle
$$ p:\E(H^{1/2}(\partial\Omega))\to\P(H^{1/2}(\partial\Omega))\,.
$$ Now, in order to write down the extensions of $A_\min$ together
with their resolvents, we make explicit the operator $G_z$ defined in
(\ref{gz}). By Theorem \ref{aggiunto}, since $A_\max=A_\min^*$, we have 
$$
\D(A_\max)=\{u=u_0+G_0h\,,\ u_0\in H^2(\Omega)
\cap H_0^1(\Omega)\,,\ h\in H^{1/2}(\partial\Omega)\}\,,
$$
$$
A_\max u=A_0u_0\,.
$$
Thus
$A_\max G_0h=0$ and so by (\ref{Green2}) there follows,
for all $h\in H^{1/2}(\partial\Omega)$ and for all $u\in\D(A_0)$,
$$ \langle G_0h,A_0u\rangle_{L^2(\Omega)}= (\hat \rho
G_0h,\tau_au)_{-\frac{1}{2},\frac{1}{2}}\,.
$$ Since, by (\ref{Green}),
\begin{align*} 
&\langle G_0h,A_0u\rangle_{L^2(\Omega)}= 
\langle G_0h,A_\max u\rangle_{L^2(\Omega)}\\=&
\langle G_0h,A_\min^*u\rangle_{L^2(\Omega)}=
\langle
h,\tau_au\rangle_{H^{1/2}(\partial\Omega)}\,,
\end{align*} 
one obtains $\hat \rho G_0h=\Lambda h$, where
$$ \Lambda:H^{1/2}(\partial\Omega)\rightarrow H^{-1/2}(\partial\Omega)
$$ is the unitary operator defined by
$$ \forall h_1,h_2\in H^{1/2}(\partial\Omega)\,,\qquad (\Lambda
h_1,h_2)_{-\frac{1}{2},\frac{1}{2}}= \langle
h_1,h_2\rangle_{H^{1/2}(\partial\Omega)}\,.
$$ For successive notational convenience we pose
$\Sigma:=\Lambda^{-1}$.
\begin{remark} If $\partial\Omega$ carries a Riemannian structure then 
$H^s(\partial \Omega)$ can be defined as the completion of
  $C^\infty(\partial\Omega)$ with respect of the scalar product
$$ \langle f,g\rangle_{H^s(\partial\Omega)} :=\langle
  f,(-\Delta_{LB}+1)^{s}g\rangle_ {L^2(\partial\Omega)}\,.
$$ Here the self-adjoint operator $\Delta_{LB}$ is the
  Laplace-Beltrami operator in $L^2(\partial\Omega)$.  With such a
  definition $(-\Delta_{LB}+1)^{1/2}$ can be extended to the unitary
  map $\Lambda$.
\end{remark}
Since $G_z=G_0+zA_0^{-1}G_z$ by (\ref{GR1}), $G_zh$ is the solution of
the Dirichlet boundary value problem
\begin{equation}\label{dbvp}
\begin{cases}
A_\max G_zh=zG_zh\,,\\ \hat\rho\, G_zh=\Lambda h\,.
\end{cases}
\end{equation}
Thus we can write $G_0\Sigma=K$, where
$K:H^{-1/2}(\partial\Omega)\to\D(A_\max)$ is the Poisson operator with
provides the solution of the Dirichlet problem with boundary data in
$H^{-1/2}(\partial\Omega)$.  Analogously we define
$K_z:H^{-1/2}(\partial\Omega)\to\D(A_\max)$ by $K_z:=G_z\Sigma$. Note
that $G_0h$, hence $G_zh$, is uniquely defined as the solution of
(\ref{dbvp}): for any other solution $u$ one has
$u-G_0h\in\K(A_0)=\{0\}$.  \par Now, according to (\ref{gamma}), we
define the bounded linear operator
\begin{equation*}
\Gamma_z: H^{1/2}(\partial\Omega)\to H^{1/2}(\partial\Omega)\,,
\quad\Gamma_z:=\tau(G_0- G_z)\,,
\end{equation*}
which, by (\ref{GR1}) and the definitions of $K$ and $K_z$, can be re-written as
\begin{equation}\label{gammaz} \begin{aligned}
\Gamma_z=-z\tau_aA_0^{-1}G_z \equiv z\,\tau_a
    R_zK\Lambda \equiv(\hat\tau_a K-\hat\tau_a K_z)\Lambda\,.
\end{aligned}\end{equation}
By $\hat\rho\, G_0h=\Lambda h$, by Theorem \ref{aggiunto} and by
Remark \ref{regtau}, we can define the regularized trace operator
$$ \hat \tau_{a,0}:\D(A_\max)\to H^{1/2}(\partial\Omega)\,,
$$
\begin{equation}\label{tau}
\hat\tau_{a,0}u:=\tau_a(u-G_0\Sigma\hat\rho u)\equiv \hat\tau_a
u-P_a\hat\rho u\equiv \tau_aA_0^{-1}A_\max u\,,
\end{equation}
where the linear operator $P_a$, known as the Dirichlet-to-Neumann
operator over $\partial\Omega$, is defined by
$$ P_a:H^{-1/2}(\partial\Omega)\to H^{-3/2}(\partial\Omega)\,,\quad
P_a:=\hat\tau_a\, K\,.
$$ In conclusion, by Theorems \ref{estensioni} and \ref{bundle}, one
has the following
\begin{theorem}\label{elliptic} Any self-adjoint extension $\hat A$ of $A_\min$ 
is of the kind
$$ \hat A:\D(\hat A)\subseteq L^2(\Omega)\to L^2(\Omega)\,,\qquad \hat
Au=A_\max u\,,
$$
\begin{align*}
\D(\hat A) =\left\{u\in \D(A_\max)\,:\,\Sigma\hat\rho
u\in\D(\Theta)\,, \quad\Pi\hat\tau_{a,0} u=\Theta \Sigma\hat\rho u
\right\}\,,
\end{align*}
where $(\Pi,\Theta)\in\E(H^{1/2}(\partial\Omega))$, and
\begin{equation*}
(-\hat A+z)^{-1} =(-A_0+z)^{-1}+
  G_z\Pi\,(\Theta+\Pi\,\Gamma_z\Pi)^{-1}\Pi G_{\bar z}^*\,,
\end{equation*}
with $\tau_{a,0}$, $G_z$ and $\Gamma_z$ defined by (\ref{tau}),
(\ref{dbvp}) and (\ref{gammaz}) respectively.
\end{theorem}
\begin{remark}
By proceding as in \cite{[Po4]}, Example 5.5, in the case the
$L^2(\partial\Omega)$-symmetric, bounded linear operator
$B:H^{3/2}(\partial\Omega)\to H^{1/2}(\partial\Omega)$ is such that
$\Theta_B:=(-P_a+B)\Lambda$, $\D(\Theta_B) =H^{5/2}(\partial\Omega)$,
is self-adjoint ($B$ pseudo-differential of order strictly less than one
suffices), the extension $A_B$ corresponding to $(\uno,\Theta_B)$ has
domain defined by Robin-type boundary conditions:
$$ \D(A_B):=\{u\in H^{2}(\Omega)\,:\, \tau_au=B\rho\}\,.
$$
\end{remark}

\end{section}

\begin{section}{A simple example.}
One of the simplest examples is given by a rotation invariant 2nd order
elliptic differential operator on the unit disc $D\subset \RE^2$. Thus
we consider the self-adjoint extensions of
$$ A_{\min}:H_0^2(D)\subset L^2(D)\rightarrow L^2(D)\,,\qquad A_\min
u=A u
$$ where
$$ A=\nabla\!\cdot \!a\nabla-c\,, \qquad
a_{ij}(x)=a\left(\|x\|\right)\, \delta_{ij}\,,\quad
c(x)=c\left(\|x\|\right)\,.
$$ We suppose that $a$ is Lipschitz continuous, $\inf_{0\le r\le 1}
a(r)>0$, and that $c\in L^q((0,1);rdr)$, $q>2$. By adding, if necessary, 
a constant to $c$ we suppose that $-A_0>0$.\par In $L^2(D)\simeq
L^2((0,1);rdr)\otimes L^2((0,2\pi);d\varphi)$ we use the orthonormal
basis $\{U_{mn}\}$, $m\in \NA$, $n\in\IN$,
$$ U_{mn}(r,\varphi) =u_{m|n|}(r)\,\frac{e^{in\varphi}}{\sqrt{2\pi}}
\,.
$$ made by the normalized eigenfunctions of the Friedrichs extension
$A_0$ of $A$. Here $\{u_{mn}\}$, $m\in \NA$, is the orthonormal basis
in $L^2((0,1);rdr)$ made by the normalized eigenfunctions of the
self-adjoint Sturm-Liouville operator
$$ L_nf(r)=-\frac{1}{r} \left(ra(r)f'(r)\right)'+
\left(c(r)+\frac{n^2}{r^2}\right)f(r)\,,\quad n\ge 0\,,
$$ with boundary conditions $f(1_-)=0$ and $\lim_{r\to 0_+} r f'(r)=0$
if $n=0$, $f(0_+)=0$ if $n\not=0$.  Denoting by $\lambda^2_{mn}>0$,
$m\in \NA$, the eigenvalues of $L_n$, one has
$$\sigma(A_0)=\sigma_d(A_0)=\{-\lambda_{m|n|}^2\,,\ m\in \NA\,,
n\in\IN\}\,.$$ 
In $H^{1/2}(S^1)$ we use the orthonormal basis $\{e_k\}$,
$k\in\IN$, defined by
$$ {e}_k(\varphi):=\frac{e^{ik\varphi}}{\sqrt{2\pi}(k^2+1)^{1/4}}\,.
$$ We want now to compute the matrix elements, relative to the basis
$\{U_{mn}\}$, of the resolvents of the self-adjoint extensions of
$A_\min$.\par By defining
$$ \nu_{mn}:= \lim_{r\uparrow 1}\, a(r)\,u_{mn}' (r)\,,
$$ one has
\begin{align*}
[G_0]_{mnk}:=&\langle U_{mn}, G_0e_k\rangle_{L^2(D)} =\langle G_0^*
U_{mn}, e_k\rangle_{H^{1/2}(S^1)} =:\overline{[G_0^*]_{kmn}}\\ =&
\langle \tau_a (-A_0)^{-1}U_{mn}, e_k\rangle
_{H^{1/2}(S^1)}=(n^2+1)^{1/4}\,
\frac{\nu_{m|n|}}{\lambda_{m|n|}^2}\ \delta_{nk}\,.
\end{align*}
Since $G_z=G_0-z(-A_0+z)^{-1}G_0$, one has then
\begin{align*}
&[G_z]_{mnk}=\overline{[G_z^*]_{kmn}}
  =[G_0]_{mnk}-\frac{z}{\lambda_{m|n|}^2+z}\,[G_0]_{mnk}\\ =&\frac{\lambda^2_{m|n|}}{\lambda^2_{m|n|}+z}\,[G_0]_{mnk}=
  (n^2+1)^{1/4}\,\frac{\nu_{m|n|}}{\lambda^2_{m|n|}+z}
  \ \delta_{nk}\,.
\end{align*}
Analogously
\begin{align*}
[\Gamma_z]_{ik}:=&-z\langle e_i,\tau_a(-A_0+z)^{-1}G_0
e_k\rangle_{H^{1/2}(S^1)}\\ =& -z(k^2+1)^{1/2}\sum_{m=1}^{\infty}
\frac{\nu_{m|k|}^2}{\lambda^2_{m|k|}(\lambda_{m|k|}^2+z)}
\,\,\delta_{ik}\,.
\end{align*}
Thus, in the case the orthogonal projection $\Pi$ is the one
corresponding to the subspace of $H^{1/2}(S^1)$ generated by
$\{e_k\,,\ k\in I\}$, $I\subseteq {\mathbb Z}$, and
$[\Theta]_{ik}=\theta_k\delta_{ik}$, $k\in I$, by Theorem
\ref{estensioni} one obtains
\begin{align*}
&\Big[(-A_{\Pi,\Theta}+z)^{-1}\Big]_{mn\tilde m\tilde n} :=
\langle U_{mn},
  (-A_{\Pi,\Theta}+z)^{-1}
  U_{\tilde m\tilde n}\rangle_{L^2(D)}\\ =&\frac{\delta_{m\tilde m}
\delta_{n\tilde n}}{\lambda^2_{m|n|}+z}
  +\frac{(n^2+1)^{1/2}}{\theta_n+[\Gamma_z]_{nn}} 
\,\,\frac{\nu_{m|n|}}
  {\lambda^2_{m|n|}+z}\,\frac{ \nu_{\tilde m|n|}}
  {\lambda^2_{\tilde m|n|}+z}\,\,\delta_{n\tilde n}
\end{align*}
for any $n\in I$, and
    \begin{equation*}
    \Big[(-A_{\Pi,\Theta}+z)^{-1}\Big]_{mn\tilde m\tilde n}
    =\frac{\delta_{m\tilde m}\delta_{n\tilde n}}{\lambda^2_{m|n|}+z}
    \end{equation*}
    for any $n\notin I$. Once the resolvent has been written as above, by 
    Theorem \ref{spettro} given any
    sequence 
    $$\{\lambda_n\}_{n\in I}\subset \mathbb{R}\cap \rho(A_0)\,,$$  
    posing $$\theta_n:=-[\Gamma_{\lambda_n}]_{nn}\,,\quad n\in I\,,$$ 
    one obtains  
    $$
    \{\lambda_n\}_{n\in I}\subset
    \sigma_p(A_{\Pi,\Theta})\,.
$$
    Moreover   
    $$
    U_n=
    \sum_{m=1}^\infty\frac{\nu_{m|n|}}{\lambda^2_{m|n|}+\lambda_n}\,
    \,U_{mn}\,,
    $$ 
    is eigenfunction with eigenvalue $\lambda_n$.

\begin{remark}\label{decomposable}
The previous example can be re-phrased in the language of decomposable 
operators (see e.g. \cite{[RS]}, section XIII.16): the operator $A_0$ is 
decomposable with fibers $A_0(n)=-L_{|n|}$ and the decomposable self-adjoint extensions 
of $A_\min$ have decomposable  
resolvents with fibers given by the resolvents of the 
self-adjoint extensions of the fibers $A_\min(n)$, which are symmetric operators 
with deficiency indices $(1,1)$. However this approach gives a 
less (than the one provided by Theorem \ref{elliptic}) explicit espression for the self-adjointness domain.
\end{remark}
\begin{remark}
In the case $a=1$, $c=0$, one has   $$\lambda_{mn}=\mu_{mn}\,,\quad
u_{mn}(r)=c_{mn}J_{n}(\mu_{mn}r)\,,$$ where  $J_n$ denotes the $n$-th
order Bessel function, $\mu_{mn}$  is its $m$-th positive zero and
$c_{mn}$ is the  normalization constant. Thus
$$\nu_{mn}=-c_{mn}\mu_{mn}J_{n+1} (\mu_{mn})\,.$$
\end{remark}
The following remark 
shows that the boundary conditions corresponding to 
couples $(\Pi,\Theta)$ of the kind above can be quite different 
from the usual ones.
\begin{remark}\label{remark} 
Suppose in the previous example we take $a=1$, $c=0$, i.e. 
$A=\Delta$, and $I=\{0\}$, 
$\lambda_0=0$. Then 
$$
\Pi: H^{1/2}(S^1)\to\CO\,,\qquad \Pi\,f=\frac{1}{\sqrt{2\pi}}\int_0^{2\pi} 
f(\varphi)\,d\varphi
$$
and $\Theta:\CO\to\CO$ is the multiplication by zero, since $\Gamma_0=0$. 
Thus
$$
\D(\Delta_{\Pi,0})
=\left\{u\in\D(\Delta_\max)\,:\,\Sigma\hat\rho u=\text{\rm const}\,,\ 
\int_0^{2\pi} 
\hat \tau_{1,0}u\,(\varphi)\,d\varphi=0\right\}\,.
$$ 
Since $\Lambda\equiv \Sigma^{-1}$ maps constants into constants, 
$$ 
\left\{u\in\D(\Delta_\max)\,:\,\hat\rho u=\text{\rm const}\right\}=
\left\{u\in H^2(D)\,:\,\rho u=\text{\rm const}\right\}
$$
by elliptic regularity, and  
$$
\int_0^{2\pi} 
\hat \tau_{1,0}u\,(\varphi)\,d\varphi=
\int_0^{2\pi}[\tau_{1}\Delta_0^{-1}\Delta_\max u](\varphi)\,d\varphi=
\int_0^{2\pi}\tau_{1}u\,(\varphi)\,d\varphi
\ \,,
$$
in conclusion one has
$$
\D(\Delta_{\Pi,0})
=\left\{u\in H^2(D)\,:\,\rho u=\text{\rm const}\,,\ 
\int_0^{2\pi} \rho\,\frac{\partial u}{\partial r}\,(\varphi)\,d\varphi=0
\right\}\,.
$$ 

\end{remark}

\end{section}


\begin{thebibliography}{99}

\bibitem{[BL]} J. Behrndt, M. Langer: 
Boundary value problems for elliptic partial differential operators on bounded domains. 
{\it J. Funct. Anal.} {\bf 243} (2007), 536-565. 

\bibitem{[BS]} M.Sh. Birman, G.Ye. Skvortsov: On the square
  summability of the highest derivatives of the solution to the
  Dirichlet problem in a region with piecewise smooth 
boundary. {\it Izv. Vyssh. Uchebn. Zaved. Matem.} {\bf 30} (1962), 12-21 
[In Russian] 

\bibitem{[BGW]} B.M. Brown, G. Grubb, I.G. Wood: M-functions for closed extensions of adjoint pairs of operators with applications to elliptic boundary problems. Preprint 2008, arXiv:0803.3630


\bibitem{[BMNW]} B.M. Brown, M. Marletta, S. Naboko, I.G. Wood: Boundary triplets and M-functions for non-selfadjoint operators, with applications to elliptic PDEs and block operator matrices. arXiv:0704.2562, to appear in 
{\it J. London Math. Soc.}

\bibitem{[C]} J.W. Calkin: General Self-Adjoint Boundary Conditions for Certain Partial Differential Operators. {\it Proc. Nat. Acad. Sci. U.S.A.} {\bf 25} (1939), 201-206
   

\bibitem{[DM]} V.A. Derkach, M.M. Malamud: Generalized Resolvents and
the Boundary Value Problem for Hermitian Operators with Gaps. 
{\it J. Funct. Anal.} {\bf 95} (1991), 1-95 

\bibitem{[EE]} D.E. Edmund, W.D. Evans: {\it Spectral 
Theory and Differential Operators.} Oxford Univ. Press, 1987

\bibitem{[GM1]} F. Gesztesy, M. Mitrea: Robin-to-Robin Maps and Krein-Type Resolvent Formulas for Schr\"odinger Operators on Bounded Lipschitz Domains. Preprint 2008, arXiv:0803.3072

\bibitem{[GM2]} F. Gesztesy, M. Mitrea: Generalized Robin Boundary Conditions, Robin-to-Dirichlet Maps, and Krein-Type Resolvent Formulas for Schr\"odinger Operators on Bounded Lipschitz Domains. Preprint 2008, arXiv:0803.3179

\bibitem{[GG]} V.I. Gorbachuk, M.L. Gorbachuk: {\it Boundary Value Problems for Operator Differential Equations.} Kluwer Academic, 1991. 

\bibitem{[Gr]} P. Grisvard: \emph{Elliptic Problems in Nonsmooth Domains}. Pitman, 1985.

\bibitem{[G]} G. Grubb:  A characterization of the non local boundary value problems associated with an elliptic operator. \emph{Ann. Scuola Norm. Sup. Pisa Cl. Sci.} \textbf{22} (1968), 425-513

\bibitem{[LU]} O.A. Ladyzhenskaya, N. N. Ural'tseva: \emph{Linear and Quasilinear Elliptic Equations}. Academic Press, 1968.

\bibitem{[LM]} J.L. Lions, E. Magenes: \emph{Non Homogeneous Boundary Value Problems and Applications, vol. I}. Springer-Verlag, 1972.

\bibitem{[vN]} J. von Neumann: Allgemeine Eigenwerttheorie Hermitscher Funktionaloperatoren. {\it Math. Ann.} {\bf 102} (1929-30), 49-131

\bibitem{[Po1]} A. Posilicano: A Kre\u\i n-like Formula for Singular
Perturbations of Self-Adjoint Operators and Applications. 
{\it J. Funct. Anal.} {\bf 183} (2001), 109-147

\bibitem{[Po2]} A. Posilicano: Self-Adjoint Extensions by Additive
Perturbations. {\it Ann. Scuola Norm. Sup. Pisa Cl. Sci.}(5) {\bf 2}
(2003), 1-20

\bibitem{[Po3]} A. Posilicano: Boundary Triples and Weyl Functions for
  Singular Perturbations of Self-Adjoint Operators. {\it
  Methods. Funct. Anal. Topology} {\bf 10} (2004), 57-63 

\bibitem{[Po4]} A. Posilicano: Self-Adjoint Extensions of
  Restrictions. arXiv:math-ph/0703078, to appear in \emph{Operators and Matrices} 
\textbf{2}
(2008).

\bibitem{[Ps]} O. Post: First order operators and boundary triples. 
{\it Russ. J. Math. Phys.}  {\bf 14}  
(2007), 482-492

\bibitem{[RS]} M. Reed, B. Simon: \emph{Methods of Modern Mathematical Physics IV. Analysis of Operators}. Academic Press, 1978.

\bibitem{[Ry]} V. Ryzhov:
 A general boundary value problem and its Weyl function.
	 {\it Opuscula Math.} {\bf 27} (2007), 305-331       

\bibitem{[V]} M.I. Vishik: On general boundary value problems for elliptic differential equations. {\it Amer. Math. Soc. Transl.} 
{\bf 24} (1963), 107--172

\end{thebibliography}
\end{document}